\documentclass[12pt,a4paper]{article}
\usepackage[T1]{fontenc}
\usepackage[utf8]{inputenc}
\usepackage{amsmath,amssymb,amsthm,mathtools}
\usepackage[margin=3cm]{geometry}
\usepackage{hyperref}

\theoremstyle{plain}
\newtheorem{theorem}{Theorem}[section]
\newtheorem{proposition}[theorem]{Proposition}
\newtheorem{lemma}[theorem]{Lemma}

\theoremstyle{definition}
\newtheorem{definition}[theorem]{Definition}
\newtheorem{remark}[theorem]{Remark}

\newtheorem*{maintheorem}{Theorem~\ref*{thm:main}}
\newtheorem*{structuretheorem}{Theorem~\ref*{thm:structure}}

\DeclareMathOperator{\tr}{tr}

\newcommand{\U}{\mathcal{U}}
\newcommand{\V}{\mathcal{V}}
\newcommand{\N}{\mathcal{N}}
\newcommand{\LL}{\mathcal{L}}
\newcommand{\RR}{\mathrm{R}}
\newcommand{\KK}{\mathcal{K}}
\newcommand{\nab}{\nabla}

\begin{document}
\title{On Riemannian $r$-th order recurrent spaces}

\author{Jos\'e M. M. Senovilla$^{1,2,3}$\\
$^1$Departamento de F\'isica, Universidad del Pa\'is Vasco UPV/EHU,\\
Apartado 644, 48080 Bilbao, Spain\\
$^2$ EHU Quantum Center, Universidad del Pa\'{\i}s Vasco UPV/EHU.\\
$^3$ e-mail: josemm.senovilla@ehu.eus}

{\let\newpage\relax\maketitle}

\begin{abstract}
We prove that a connected Riemannian manifold of dimension $n\ge3$ with
irreducible local holonomy which is $r$-th order recurrent,
$\nab^r\RR=t^{(r)}\otimes\RR$, is locally symmetric, for every $r\ge1$. Without
any assumption on the holonomy, we show that such a manifold is either locally
symmetric or, locally, a flat extension $\Sigma^2\times\mathbb R^{n-2}$ of a
surface: the only possible exceptions to local symmetry are essentially two-dimensional.
\end{abstract}

\maketitle

\section{Introduction}

Let $(M,g)$ be an $n$-dimensional ($n\geq 3$) connected semi-Riemannian manifold with curvature tensor $\RR$ and
Levi-Civita connection $\nab$. It is called \emph{$r$-th order recurrent} if
\begin{equation}\label{eq:def}
\nab^r \RR =t^{(r)} \otimes \RR
\end{equation}
for some $r$-covariant tensor field $t^{(r)}$. For $r=1$ this is the
classical notion of a recurrent space \cite{Ru,walker2}, while $t^{(r)}=0$
gives the $r$-th order symmetric spaces \cite{senovilla2008,kaigorodov}; no
condition on $t^{(r)}$ is imposed, so that the latter are included. The word
\emph{proper} will be used in the sense of \cite[Definition
4]{senovilla}: an $r$-th order recurrent $(M,g)$ is proper of order $r$ if $r$ is the least order such that \eqref{eq:def} holds.

The condition \eqref{eq:def} has a long history. The $r$-th order symmetric
case $\nab^r\RR=0$ was settled in the Riemannian signature by Lichnerowicz
\cite{lichnerowicz}, Nomizu and Ozeki \cite{nomizu} and Tanno \cite{tanno},
proving that only local symmetry is permitted; the Lorentzian second-order
symmetric case was opened in \cite{senovilla2008}. For $r=1$, Riemannian
recurrent spaces of dimension $n>2$ are locally symmetric, apart from the flat
extensions of surfaces discussed below \cite{mogi,walker2}, while the
semi-Riemannian classification is due to Walker \cite{walker2}; proper
second-order recurrent Riemannian spaces were studied by Roter
\cite{roter1,roter2}. The case
$r=2$ has been studied in detail: Thompson proved that a proper second-order
recurrent manifold which is conformally flat is actually recurrent
\cite{thompson68}, and that the same conclusion holds when the scalar
curvature does not vanish and either the Ricci tensor is definite or $n=3$
(and also for $n=4$ Lorentzian) \cite{thompson70}; the four-dimensional
Lorentzian case was completed by McLenaghan and Thompson
\cite{mclenaghan}, and second-order conformally recurrent plane-fronted waves
were examined by Collinson and S\"oler \cite{collinson}. These results are
collected in Proposition 12 of \cite{senovilla}. 

Spaces satisfying \eqref{eq:def} form a distinguished subclass of the
\emph{$KS_n$ spaces} introduced and investigated  in \cite{senovilla,kaigorodov}, a considerably more
general family for which the structure theory used below is developed. All results we import are stated
there for $KS_n$ spaces and therefore apply a fortiori to the $r$-th
order recurrent case treated here. In particular, proper $r$-recurrent spaces do exist for arbitrary signature and explicit examples are exhibited in \cite{senovilla} for all $r\geq 1$ and for all dimensions $n> 2$, {\em except for the case of positive definite metric}. In this Riemannian signature, Theorem 6 of \cite{senovilla} (due to
\cite{lichnerowicz,tanno} and \cite{roter1,roter2}) states that recurrent (up to
flat extensions of surfaces) and proper second-order recurrent manifolds with
$n>2$ are locally symmetric, and
a conjecture is put forward there, for the whole $KS_n$ class, to the effect
that the same should hold in general. 

In this work, I address this conjecture in the outstanding case
\eqref{eq:def}. The conjecture requires dimension $n>2$, because every surface
with nowhere vanishing non-constant Gaussian curvature is $r$-recurrent. Notice, though, that they are {\em proper} recurrent of order $r=1$, and so are its flat extensions $\Sigma^2\times\mathbb R^{n-2}$, in every
dimension $n>2$; none of them is locally symmetric. However, flat extensions of surfaces may also be proper of higher
order if the Gaussian curvature can vanish: Remark \ref{rem:sharp} exhibits one which is proper of order two.
Hence the conjecture can only hold under an assumption that excludes the flat
extensions of surfaces. Since, by de Rham's theorem, every Riemannian manifold
is locally a product of a flat factor and irreducible factors, the natural
assumption is the irreducibility of the local holonomy. Our main result
confirms the conjecture in this form.

\begin{maintheorem}
Let $(M,g)$ be a connected Riemannian manifold of dimension $n\ge3$ whose local
holonomy group acts irreducibly at every point. If $(M,g)$ is $r$-th order
recurrent for some $r\ge1$, then it is locally symmetric.
\end{maintheorem}

Without any assumption on the holonomy, flat extensions of two-dimensional surfaces $\Sigma^2$ are the only
possible exceptions. Accordingly, we shall say that $(M,g)$ is \emph{locally a
flat extension of a surface} if its curvature operator has rank at most one and
every point where $\RR\neq0$ has a neighbourhood isometric to a Riemannian
product $\Sigma^2\times\mathbb R^{n-2}$. We prove:

\begin{structuretheorem}
Let $(M,g)$ be a connected Riemannian manifold of dimension $n\ge3$ which is
$r$-th order recurrent, $\nab^{r}\RR=t\otimes\RR$, for some $r\ge1$. Then
$(M,g)$ is either locally symmetric or locally a flat extension of a surface.
More precisely, if $(M,g)$ is not locally symmetric, then:
\begin{enumerate}
\item[(a)] on $\operatorname{int}\{t=0\}$ one has $\nab\RR=0$;
\item[(b)] the curvature operator has rank at most one on all of $M$;
\item[(c)] on the open set $\{\RR\neq0\}$, $(M,g)$ is locally isometric to a
  product $\Sigma^2\times\mathbb{R}^{n-2}$. On $\{t\neq0\}\cap\{\RR\neq0\}$ the
  curvature $K$ of $\Sigma^2$ satisfies $\nab^{r}K\neq0$, whereas on the
  components of $\operatorname{int}\{t=0\}$ where $\RR\neq0$, $K$ is constant.
\end{enumerate}
\end{structuretheorem}

In Riemannian signature, $r$-th order recurrence beyond local symmetry is thus a
two-dimensional phenomenon. Theorem \ref{thm:structure} is considerably stronger
than what can be obtained from the general theory of $KS_n$ spaces. Theorem 1 of
\cite{senovilla} shows that a decomposable $KS_n$ space with $t^{(r)}\neq0$ can
only be a flat extension, but it says nothing about the indecomposable ones, and
it allows the non-flat factor to be of any dimension. Theorem
\ref{thm:structure} shows that, for $r$-th order recurrent Riemannian manifolds,
no decomposability assumption is needed, and that the non-flat factor is
necessarily \emph{two-dimensional}: the only exceptions to local symmetry are the
(local) flat extensions of surfaces, a phenomenon already present in dimension
two. The result is sharp, since such exceptions exist (Remark
\ref{rem:sharp}).

The proof reduces the problem to conullity two with
nilpotent splitting tensor, and then to a statement about a function of one
variable: along the integral curves of a privileged direction, certain
components of $\nab^r\RR$ are governed by an equation of type $y''+Ky=0$, where $K$ is (up to a constant) the scalar curvature of the manifold, and then $r$-th order recurrence forces a rational solution whose poles are incompatible with the
recurrence itself.

\section{Reduction to conullity two}

For $(M,g)$ satisfying \eqref{eq:def} with some $r\ge1$, we set
\[
 \U:=\{p\in M:\ t^{(r)}_p\neq0\ \text{and}\ \RR_p\neq0\}.
\]
\begin{proposition}[\cite{senovilla}]
\label{prop:reduction}
Let $(M,g)$ be a Riemannian manifold satisfying \eqref{eq:def} with $r\ge2$.
Then, on $\U$, the curvature tensor has the form
\begin{equation}\label{eq:typeII}
\RR=K\,F\otimes F,
\end{equation}
where $F$ is a simple unit $2$-form and $K$ is a
nowhere-vanishing function proportional to the scalar
curvature.
\end{proposition}

\begin{proof}
By Proposition 2 and Theorem 2 of \cite{senovilla},
at points where $t^{(r)}\neq0$ and $\RR\neq0$, $\dim\mathcal L$ is either
one or two, where $\mathcal L$ is the module of
$1$-forms $\omega$ satisfying
\[
\omega\wedge \RR=0.
\]
Proposition 8 and Corollary 6 of \cite{senovilla}
exclude type $\mathrm I_\epsilon$ for $r\geq2$.
Moreover, type $\mathrm I_N$ is impossible in
Riemannian signature. Hence only type $\mathrm{II}_1$
is possible; it has $\dim\LL=2$, and the result follows from Lemma 5
of \cite{senovilla}.
\end{proof}

From now on we work on an open set where \eqref{eq:typeII} holds with $K\neq0$ (i.e. $\RR \neq 0$).
Write
\[
 \N := \LL^{\perp} = \{X: \RR(X,\cdot,\cdot,\cdot)=0\},
\]
so $\dim\N=n-2$: the manifold has \emph{conullity two}. We use the sign
convention
\[
 \RR(X,Y,Z,W)=g(\RR(X,Y)W,Z),\qquad \RR(X,Y)=[\nab_X,\nab_Y]-\nab_{[X,Y]},
\]
so that $\RR(X,Y,X,Y)$ is the
sectional curvature of an orthonormal pair $X,Y$. Since
$(F\otimes F)(e_1,e_2,e_1,e_2)=1$ for any orthonormal basis $e_1,e_2$ of $\LL$
with $F=e^1\wedge e^2$, $K$ is the sectional curvature of the plane field
$\LL$, and the scalar curvature is $2K$. (In \cite{VanHook} the same quantity
is written $R(e_1,e_2,e_2,e_1)=\frac12\mathrm{Scal}$.) Riemannian manifolds of
conullity two have been studied systematically in \cite{bkv}; they also form
one of the classes in Szab\'o's structure theory of semi-symmetric spaces
\cite{szabo}, which is relevant here because type $\mathrm{II}$ spaces are
semi-symmetric (Theorem 4 of \cite{senovilla}).

\begin{lemma}\label{lem:parallelL}
$\N$ is an involutive distribution whose leaves are flat and
totally geodesic. Along any leaf, $\LL$ is a parallel subbundle
and its unit area form $F$ is parallel.
\end{lemma}

\begin{proof}
That the nullity distribution is involutive with totally
geodesic leaves, on the open set where its dimension is
constant, is classical \cite{ChernKuiper,Maltz}. The leaves
are flat by the Gauss equation, since they are totally
geodesic and $\RR(U,\cdot,\cdot,\cdot)=0$ for $U\in\N$.

For $\ell\in\Gamma(\LL)$ and $U,V\in\Gamma(\N)$, we have
\[
g(\nabla_U \ell,V)
=U(g(\ell,V))-g(\ell,\nabla_U V)=0,
\]
since $\N$ is totally geodesic. Hence
$\nabla_U \ell\in\Gamma(\LL)$, so $\LL$ is parallel along each
leaf of $\N$.

The induced connection on $\LL$ along each leaf is metric.
Since $\dim\LL=2$, it preserves the unit area form $F$.
\end{proof}

\begin{definition}
The \emph{splitting tensor} $C:\N\times\LL\to\LL$ is
$C(U,X):= C_U X := -(\nab_X U)^{\LL}$, where the super-script $\LL$ stands for projection to $\LL$.
\end{definition}

$C$ is tensorial in both arguments, and $C\equiv0$ if and only if $\N$ (hence
also $\LL$) is parallel, i.e.\ if and only if $(M,g)$ is locally decomposable.

Two classical identities will be used, see \cite{VanHook,VanHookThesis}. Along a unit-speed
geodesic $\gamma(s)$ tangent to $\N$, let $U=\dot\gamma$
and identify the spaces $\LL_{\gamma(s)}$ by parallel
transport along $\gamma$. Then
\begin{equation}\label{eq:riccati}
 C(s)=C_0(I-sC_0)^{-1},
 \qquad
 K(s)=\frac{K_0}{\det(I-sC_0)} ,
\end{equation}
where $K$ is the sectional curvature of $\LL$, $K_0=K(0)$ and $C_0=C(0)$.
The first identity follows from a Riccati equation for $C$, while the second follows from the
second Bianchi identity, which gives
$U(K)=K\tr C_U$ \cite{VanHook,VanHookThesis}.

\begin{lemma}[Nilpotency]\label{lem:nilpotent}
If $r\geq 2$, for every $U\in\N$ the endomorphism $C_U$ of $\LL$ is nilpotent, and $K$ is
constant along the leaves of $\N$.
\end{lemma}

\begin{proof}
By Lemma 3 and Corollary 4 of \cite{senovilla}, since
$t^{(r-1)}=0$, the tensor $t^{(r)}$ is symmetric in its
last two indices and its last index annihilates $\N$.
Consequently,  the defining recurrence equation \eqref{eq:def} gives in particular
\[
\nabla^r \RR (U,\cdots , U;\cdot,\cdot,\cdot,\cdot) =0
\qquad\text{for every }U\in\N.
\]

Let $\gamma$ be a unit nullity geodesic, $U=\gamma'$. Since
$\nab_{\gamma'}\gamma'=0$, contracting all $r$ derivative slots with $\gamma'$
gives the $r$-th iterated derivative of $\RR$ along $\gamma$; in a parallel frame along
$\gamma$ this is $d^r/ds^r$ of the frame components. By Lemma
\ref{lem:parallelL}, $F$ is parallel along $\gamma$, so by \eqref{eq:typeII}
the frame components of $\RR$ are $K(s)$ times constants. Hence
$$\frac{d^r K(s)}{ds^r} \equiv 0$$ 
so that \emph{$K$ is a polynomial in $s$ of degree $\le r-1$.}

By \eqref{eq:riccati}, $K(s)\,p(s)=K_0$ on a neighbourhood of $s=0$, where
$p(s)=\det(I-sC_0)$ is a polynomial of degree $\le2$. Since the manifold is of type II and is non-flat, $K_0\ne0$.
Thus $K$ and $p$ are nonzero polynomials. Since
$K(s)p(s)=K_0$ on an interval, both must be constant.
Hence $p\equiv1$ and $K$ is constant along $\gamma$. Now $p\equiv1$ means
$\tr C_0=\det C_0=0$, i.e.\ $C_0$ is nilpotent by Cayley--Hamilton, since \(\dim\LL=2\).
\end{proof}

\begin{remark}
The last step is purely local: no completeness of the nullity geodesics is
needed, because the identity between polynomials propagates from an arbitrarily
short interval. This also rules out metric cones directly: for
$g=dt^2+t^2h$ one has $C_{\partial_t}=-I/t$ and $K\propto t^{-2}$, which is not
polynomial. (This can be deduced directly and more generally from Proposition 1 of \cite{senovilla}.)
\end{remark}

\section{The adapted frame}\label{sec:frame}

Throughout this section $r\ge2$ and we work on the open set $\U$ of Proposition
\ref{prop:reduction}, so that Lemma \ref{lem:nilpotent} applies. Let
\[
 \V:=\{p\in\U:\ C_p\neq0\},
\]
an open set, possibly empty; if $\V=\emptyset$, then $C\equiv0$ on $\U$. Fix $p\in\U$. Each $C_U$, $U\in\N_p$, is a nilpotent endomorphism of the plane
$\LL_p$, and $U\mapsto C_U$ is linear, so $\{C_U:\ U\in\N_p\}$ is a linear
subspace of $\mathrm{End}(\LL_p)$ consisting of nilpotent endomorphisms. If its
dimension were at least $2$, it would meet the three-dimensional space of
self-adjoint endomorphisms of the Euclidean plane $\LL_p$ (inside the
four-dimensional $\mathrm{End}(\LL_p)$) in a non-zero element, which is
impossible because a self-adjoint nilpotent endomorphism vanishes. Hence this
subspace has dimension at most one \cite[Section 1]{VanHook}. Therefore, on $\V$,
\[
 C_U=\mu(U)\,\Theta
\]
for a nowhere vanishing $1$-form $\mu$ on $\N$ and a rank one nilpotent
$\Theta\in\mathrm{End}(\LL)$.

Locally on $\V$ we choose a unit vector field
$e_1$ spanning the line field $\ker\Theta=\operatorname{im}\Theta\subset\LL$
(it is unique up to sign), and we complete it with $e_2\in\LL$ to an orthonormal
basis such that $F=e^1\wedge e^2$. We normalise $\Theta$ by $\Theta e_2=e_1$,
which fixes $\mu$ through $C_Ue_2=\mu(U)e_1$, and we write
\[
 N:=|\mu|>0,\qquad E:=\mu^\sharp/N\in\Gamma(\N),
\]
where $\mu^\sharp\in\N$ is the metric dual of $\mu$ within $\N$.

\begin{lemma}[Adapted frame]
\label{lem:structure}
On $\V$, with the notation above and $\kappa:=g(\nab_{e_2}e_2,e_1)$:
\begin{enumerate}
\item[(i)] $\nab_{e_1}e_1=\nab_{e_1}e_2=\nab_{e_1}E=0$;
\item[(ii)] $\nab_Ue_1=\nab_Ue_2=\nab_UE=0$ for every $U\in\Gamma(\N)$;
\item[(iii)] $\nab_{e_2}e_1=-\kappa e_2+NE$, \ $\nab_{e_2}e_2=\kappa e_1$,
  \ and $\nab_{e_2}E+Ne_1\in\Gamma(\N)$;
\item[(iv)] $e_1(\kappa)=K+\kappa^2$, \ $e_1(N)=\kappa N$, \ and
  $U(\kappa)=U(N)=U(K)=0$ for every $U\in\Gamma(\N)$.
\end{enumerate}
In particular $[E,e_1]=\nab_{E}e_1-\nab_{e_1}E=0$.
\end{lemma}

\begin{proof}
Items (i), (ii) and the first two relations in (iv) are
\cite[Proposition 1.2]{VanHook} (see also
\cite[Proposition 2.1.2]{VanHookThesis}), where $T_1$, $a_1$ and $\beta$ stand
for our $E$, $N$ and $\kappa$. That proposition is stated on the set where
$C\neq0$, i.e.\ on $\V$, and its proof is a local computation that uses
neither completeness nor any further frame. The relation $U(K)=0$ is Lemma
\ref{lem:nilpotent}. For (iii), using $C_Ue_1=0$ and $C_Ue_2=\mu(U)e_1$: the
$\N$-component of $\nab_{e_2}e_1$ is determined by
$g(\nab_{e_2}e_1,U)=-g(e_1,\nab_{e_2}U)=g(e_1,C_Ue_2)=\mu(U)$, so it is
$\mu^\sharp=NE$, while its $e_2$-component is $-g(e_1,\nab_{e_2}e_2)=-\kappa$;
similarly $g(\nab_{e_2}e_2,U)=g(e_2,C_Ue_2)=0$, so
$\nab_{e_2}e_2=\kappa e_1$. Finally,
$g(\nab_{e_2}E,e_1)=-g(E,\nab_{e_2}e_1)=-N$ and
$g(\nab_{e_2}E,e_2)=-g(E,\nab_{e_2}e_2)=0$.
\end{proof}

\begin{remark}\label{rem:frenet}
On the open subset of $\V$ where the corresponding Frenet construction along
the integral curves of $e_2$ is defined, Van Hook completes $E_1:=E$ to an
orthonormal frame $\{E_1,\dots,E_{n-2}\}$ of $\N$, with functions $f_a>0$,
$f_1=N$, such that, setting $E_0:=e_1$ and $f_{n-1}:=0$,
\[
 \nab_{e_2}E_a=-f_aE_{a-1}+f_{a+1}E_{a+1},\qquad
 \nab_{e_1}E_a=\nab_{E_b}E_a=0,
\]
\[
 e_1(f_a)=\kappa f_a,\qquad E_b(f_a)=0,
\]
see \cite[Proposition 2.2]{VanHook} and
\cite[Proposition 2.2.2]{VanHookThesis}. This subset may be smaller than $\V$:
the construction breaks down, for instance, at points of $\V$ where some $f_a$
with $a\ge2$ vanishes. We shall not need this completion: all
the arguments below use only $e_1$, $e_2$, $E$ and the relations of Lemma
\ref{lem:structure}, which hold on the whole of $\V$.
\end{remark}

\begin{lemma}\label{lem:gradF}
The unit area form $F$ of $\LL$ is closed. Moreover, for
$\ell \in\LL$, $U\in\N$ and any vector field $X$,
\[
(\nabla_XF)(\ell,U)=F(\ell,C_UX),
\]
whereas $\nabla F$ vanishes when its last two arguments
belong both to $\LL$ or both to $\N$. Equivalently,
\[
\nabla F=-e^2\otimes(e^2\wedge\hat \mu),
\]
where $\hat\mu$ is the extension of $\mu$ to act on arbitrary tensor fields that vanishes on
$\LL$ (this last formula holds on $\V$). In particular, $\nabla F=0$ if and only if $C=0$.
\end{lemma}

\begin{proof}
The first statements hold wherever \eqref{eq:typeII}
holds. Since $F$ is the unit area form of $\LL$, the induced
metric connection on $\LL$ preserves $F$. Hence
\[
(\nabla_XF)(\ell_1,\ell_2)=0,
\qquad \ell_1,\ell_2\in\LL.
\]
Likewise, since $F$ vanishes whenever one of its arguments
belongs to $\N$,
\[
(\nabla_XF)(U,V)=0,
\qquad U,V\in\N.
\]

For $\ell\in\LL$ and $U\in\N$, using the definition of the
splitting tensor, we obtain
\[
(\nabla_XF)(\ell,U)
=-F(\ell,\nabla_XU)\\
=-F\big(\ell,(\nabla_XU)^{\LL}\big)\\
=F(\ell,C_UX).
\]
On $\V$ we have
$C_Ue_1=0$ and $C_Ue_2=\mu(U)e_1$, and therefore
\[
\nabla F=-e^2\otimes(e^2\wedge\mu).
\]
Consequently, antisymmetrising
\[
dF=0.
\]
Finally, $\nabla F=0$ if and only if $\mu=0$, which is
equivalent to $C=0$.
\end{proof}

\begin{proposition}\label{prop:vanishing}
On the open set where \eqref{eq:typeII} holds, $r$-th order
recurrence implies
\[
 (\nabla^{r}\RR)(X_1,\dots,X_r;\,Y_1,Y_2,Y_3,Y_4)=0
\]
whenever one of $Y_1,Y_2,Y_3,Y_4$ belongs to $\N$.

Equivalently, for every fixed choice of
$X_1,\dots,X_r$, the resulting algebraic curvature tensor
is pointwise proportional to $F\otimes F$.
\end{proposition}
\begin{proof}
Contract \eqref{eq:def} on a curvature slot with $U\in\N$: the right-hand side
is 
\[
 t^{(r)}(X_1,\dots,X_r)\,\RR(\dots,U,\dots)=0 .
\]
For the second formulation,
an algebraic curvature tensor on $\LL\oplus\N$ vanishing whenever a slot lies
in $\N$ is supported on $\LL$, and since $\dim\LL=2$ it is a multiple of
$F\otimes F$.
\end{proof}

\begin{remark}\label{rem:onlyr}
Proposition \ref{prop:vanishing} constrains $\nabla^r\RR$ only;
it gives no direct vanishing condition on $\nabla^k\RR$ for
$k<r$. In particular, the non-vanishing of
\[
(\nabla_{e_2}\RR)(e_1,e_2,e_2,E)=-KN
\]
is not in itself an obstruction when $r\geq2$.
This is why the argument must genuinely take the recurrence
order $r$ into account.
\end{remark}

\begin{remark}\label{rem:reduction}
By the recurrence equation and Proposition \ref{prop:vanishing},
we have
\[
(\nabla^r\RR)(X_1,\dots,X_r;\,\cdot,\cdot,\cdot,\cdot)
=
K\,t^{(r)}(X_1,\dots,X_r)\,F\otimes F.
\]
Thus, the curvature derivative $\nabla^r\RR$ is entirely
determined by the covariant tensor $K\,t^{(r)}$: its four
curvature slots are supported on the two-dimensional
distribution $\LL$, while all the remaining information is
encoded in $K\,t^{(r)}$.
\end{remark}

\subsection{Recovery of the first- and second-order cases}
\label{subsec:r12}

We now recover the known results for $r=1$ and $r=2$.
The first-order case is recalled for completeness, whereas
the second-order case follows directly from the geometric
framework developed above. The computation behind the second-order case
is isolated in the following lemma, which will also be used in
Section \ref{sec:proof}.

\begin{lemma}\label{lem:r2}
Let $(M,g)$ be a Riemannian manifold satisfying \eqref{eq:def} with $r=2$.
Then $C\equiv0$ on $\U$.
\end{lemma}

\begin{proof}
Suppose that the open set $\V=\{p\in\U:\ C_p\neq0\}$ is non-empty. By
Proposition \ref{prop:reduction}, $\RR=K\,F\otimes F$ with $K\neq0$ on $\U$, and
on $\V$ we use the frame of Lemma \ref{lem:structure}.

Differentiating $\RR=K\,F\otimes F$ twice, every term of
$\nabla^2\RR$ containing an undifferentiated factor $F$
vanishes on the curvature arguments $(e_2,U,e_2,U)$, with
$U\in\N$, because $F(e_2,U)=0$. The only remaining
term is $2K\,\nabla F\otimes\nabla F$, and by Lemma
\ref{lem:gradF}, $(\nabla_{e_2}F)(e_2,U)=-\mu(U)$. Hence
\[
(\nabla^2\RR)(e_2,e_2;\,e_2,U,e_2,U)
=2K\,\mu(U)^2
\qquad U\in\N.
\]
On the other hand, Proposition \ref{prop:vanishing}
implies that the left-hand side vanishes, since one of
the curvature arguments belongs to $\N$. Hence
\[
K\,\mu(U)^2=0
\qquad\text{for all }U\in\N.
\]
Since $K\neq0$ this yields $\mu=0$ on $\V$, a
contradiction. Hence $\V=\emptyset$, i.e.\ $C\equiv0$ on $\U$.
\end{proof}

\begin{theorem}\label{thm:r12}
Let $(M,g)$ be a Riemannian manifold of dimension $n\geq3$
whose local holonomy representation is irreducible at every
point. If $(M,g)$ is recurrent ($r=1$) or $2$-nd order
recurrent ($r=2$), then it is locally symmetric.
\end{theorem}

\begin{proof}
We distinguish the two cases.

\smallskip
\noindent
\textit{Case $r=1$.}
If $t^{(1)}\equiv0$, then $\nabla\RR=0$ and the manifold is
locally symmetric. That recurrent Riemannian manifolds
with irreducible local holonomy in dimension $n>2$ are
locally symmetric is classical \cite{mogi,walker2}; see also
Theorem 6(b) of \cite{senovilla}. It also follows from
Proposition \ref{prop:Czero} below.

\smallskip
\noindent
\textit{Case $r=2$.}
Let $\U$ be as before, with $r=2$. If $\U=\emptyset$, then
$t^{(2)}\otimes\RR\equiv0$, i.e.\ $\nabla^2\RR=0$, and the manifold is locally
symmetric \cite{nomizu,tanno}, see Theorem 6(a) of \cite{senovilla}. If
$\U\neq\emptyset$, Lemma \ref{lem:r2} gives $C\equiv0$ on $\U$. Thus $\N$ and
$\LL$ are parallel orthogonal distributions on $\U$, which gives a local
Riemannian product decomposition. This contradicts the assumed irreducibility of
the local holonomy representation. Therefore $\U$ must be empty, and $(M,g)$ is
locally symmetric as above.
\end{proof}

\begin{remark}
The first-order case is included only for completeness,
as it follows from the classical result on recurrent
Riemannian manifolds with irreducible local holonomy in
dimension greater than two.

The second-order case, on the other hand, admits a direct
geometric proof in the present framework. Its essential
ingredient is the quadratic contribution
$\nabla F\otimes\nabla F$ to $\nabla^2\RR$, which forces
the splitting tensor to vanish. This suggests
that the higher-order cases should be investigated by
examining the corresponding terms involving higher
covariant derivatives of $F$.
\end{remark}

\section{The chain $G_k$ along the privileged geodesics}\label{sec:chain}

Throughout this section $r\ge2$, the open set $\V$ of Section \ref{sec:frame}
is assumed to be non-empty, and we work on $\V$ with the frame of Lemma
\ref{lem:structure}. Only Lemmas \ref{lem:structure} and \ref{lem:gradF} will be used. By Lemma
\ref{lem:structure}(i),(ii), for $X\in\{e_1,E\}$ we have
$\nab_Xe_1=\nab_Xe_2=\nab_XE=0$ and $[E,e_1]=0$, and $\nab_XF=0$ by Lemma
\ref{lem:gradF}, since $\nab F=-e^2\otimes(e^2\wedge\hat\mu)$ and $e^2(X)=0$.
Hence, for any rank-$q$ tensor field $T$ evaluated on vectors
$Y_i\in\{e_1,e_2,E\}$ ($i=1,\dots,q$),
\[
 (\nab_XT)(Y_1,\dots,Y_q)=X\big(T(Y_1,\dots,Y_q)\big).
\]
Since $e_1$ is a unit vector field with $\nab_{e_1}e_1=0$, its integral curves
are unit speed geodesics. We fix one of them, defined on an interval $I$,
denote by $u$ its arclength parameter, and write ${}'=d/du=e_1$ for functions
restricted to it.

\begin{lemma}\label{lem:slots}
Let $Z_1,\dots,Z_k\in\{e_1,E\}$. Then
\[
 (\nab^{k}\RR)(Z_1,\dots,Z_k;\,\cdot\,,\cdot\,,\cdot\,,\cdot\,)
 =\big(Z_1\cdots Z_k(K)\big)\,F\otimes F .
\]
This equals 
$$
\frac{d^k K}{du^k} F\otimes F=: K^{(k)} F \otimes F \qquad \mbox{if} \quad Z_1=\dots=Z_k=e_1
$$
and vanishes if some $Z_i=E$. In particular, it vanishes whenever one of the four curvature
arguments lies in $\N$.
\end{lemma}

\begin{proof}
By induction on $k$, the case $k=0$ being \eqref{eq:typeII}. Since
$\nab_{Z_1}Z_j=0$ for all $j$,
\[
 (\nab^{k}\RR)(Z_1,\dots,Z_k;\,\cdot\,)
 =\nab_{Z_1}\big[(\nab^{k-1}\RR)(Z_2,\dots,Z_k;\,\cdot\,)\big]
 =\nab_{Z_1}\big[(Z_2\cdots Z_k(K))\,F\otimes F\big],
\]
and $\nab_{Z_1}F=0$. If some $Z_i=E$, then, since $[E,e_1]=0$ and $E(K)=0$,
the derivation $E$ commutes with the $e_1$'s and annihilates $K$.
\end{proof}

\begin{lemma}\label{lem:jacobi}
Along $I$, the function $y:=1/N$ satisfies $y''+Ky=0$, and $\kappa=-y'/y$.
\end{lemma}

\begin{proof}
By Lemma \ref{lem:structure}(iv), $e_1(N)=\kappa N$ and $e_1(\kappa)=K+\kappa^2$.
The first gives $y'=-\kappa y$, hence
$y''=-\kappa'y-\kappa y'=-(K+\kappa^2)y+\kappa^2y=-Ky$. These are the relations
of \cite[Proposition 1.2(c) and (1.2)]{VanHook}, written in our notation.
\end{proof}

Set $\mathbf W:=(e_1,e_2,E,e_2)$ and, for $k\ge0$, define along $I$
\[
 G_k:=(\nab^{k+1}\RR)(e_2,\underbrace{e_1,\dots,e_1}_{k};\;\mathbf W).
\]
Thus $G_0=(\nab_{e_2}\RR)(e_1,e_2,E,e_2)=NK$; cf.\ Remark \ref{rem:onlyr}.

\begin{lemma}\label{lem:words}
For every $k\ge1$ and $1\le i\le k$,
\[
 (\nab^{k}\RR)(\underbrace{e_1,\dots,e_1}_{i-1},e_2,
 \underbrace{e_1,\dots,e_1}_{k-i};\;\mathbf W)=\frac{d^{i-1} G_{k-i}}{du^{i-1}} :=G_{k-i}^{(i-1)} .
\]
Consequently, if $(M,g)$ satisfies \eqref{eq:def}, then for $0\le k\le r-1$ the
function $G_k$ is a polynomial in $u$ of degree at most $r-2-k$ (so
$G_{r-1}=0$).
\end{lemma}

\begin{proof}
The first $i-1$ derivative arguments are $e_1$, and $\nab_{e_1}$ acts as $e_1$
on functions of $e_1,e_2,E$; applying this $i-1$ times gives
$e_1^{\,i-1}\big((\nab^{k-i+1}\RR)(e_2,e_1,\dots,e_1;\mathbf W)\big)
=G_{k-i}^{(i-1)}$. For the second statement take $k=r$ and $i=r-j$: the
left-hand side vanishes by Proposition \ref{prop:vanishing}, since the third
curvature argument $E$ lies in $\N$; hence $G_j^{(r-j-1)}=0$ for
$0\le j\le r-1$.
\end{proof}

\begin{lemma}\label{lem:recursion}
Put $S_k:=\sum_{j=0}^{k-1}G_j^{(k-1-j)}$ (with $S_0=0$). Then
\[
 G_k=\kappa\,S_k+N\,K^{(k)},\qquad S_{k+1}=S_k'+G_k .
\]
\end{lemma}

\begin{proof}
Let $T=\nab^k\RR$ and write $\mathbf W =(e_1,e_2,E,e_2)=(W_1,W_2,W_3,W_4)$. By the Leibniz
rule,
\[
\begin{aligned}
 G_k={}&e_2\big[T(e_1,\dots,e_1;\mathbf W)\big]
 -\sum_{m=1}^{k}T(e_1,\dots,\underbrace{\nab_{e_2}e_1}_{m},\dots,e_1;\mathbf W)\\
 &-\sum_{s=1}^{4}T(e_1,\dots,e_1;W_1,\dots,\nab_{e_2}W_s,\dots,W_4).
\end{aligned}
\]
We examine each term, using Lemma \ref{lem:structure}(iii) and Lemma \ref{lem:slots}.
\begin{itemize}
\item The first term vanishes, because $T(e_1,\dots,e_1;\mathbf W)=0$, as
  $W_3=E\in\N$.
\item Derivative arguments: $\nab_{e_2}e_1=-\kappa e_2+NE$. The part $NE$ gives a
  component with $E$ in a derivative argument, which vanishes. The part
  $-\kappa e_2$ gives $-\kappa$ times the component with $e_2$ in position
  $m$, which equals $G_{k-m}^{(m-1)}$ by Lemma \ref{lem:words}. With the
  overall minus sign, the total contribution is
  $\kappa\sum_mG_{k-m}^{(m-1)}=\kappa S_k$.
\item Curvature arguments. For $s=1$, $\nab_{e_2}e_1=-\kappa e_2+NE$ gives the
  arguments $(e_2,e_2,E,e_2)$, which vanish by antisymmetry, and
  $(E,e_2,E,e_2)$, which vanish since $E\in\N$. For $s=2$, $\nab_{e_2}e_2=\kappa
  e_1$ gives $(e_1,e_1,E,e_2)$, zero by antisymmetry. For $s=4$ one obtains
  $(e_1,e_2,E,e_1)$, zero since $E\in\N$. For $s=3$, $\nab_{e_2}E=-Ne_1+U$ with
  $U\in\N$; the part $U$ gives zero, while $-Ne_1$ gives
  $-NK^{(k)}(F\otimes F)(e_1,e_2,e_1,e_2)=-NK^{(k)}$. With the overall minus
  sign this contributes $+NK^{(k)}$.
\end{itemize}
This proves the first identity; the second follows from the definition of
$S_k$.
\end{proof}

\begin{remark}
The terms in $\kappa$ coming from $\nab_{e_2}e_1=-\kappa e_2+NE$ do \emph{not}
disappear: they are exactly $\kappa S_k$, and they are what makes each $G_k$
depend on the previous ones. Only the part along $NE$ drops out.
\end{remark}

Let $p_j,q_j$ be the differential polynomials in $K$ defined by
\[
 p_0=1,\quad q_0=0,\qquad p_{j+1}=p_j'-Kq_j,\quad q_{j+1}=p_j+q_j' ,
\]
so that 
$$
\frac{d^j Y}{du^j} :=Y^{(j)}=p_jY+q_jY'$$ 
for every solution of $Y''+KY=0$. Giving a weight $2$ to  $K$ and weight $1$ to every derivative $d/du$, the polynomials $p_j$ are isobaric of weight $j$, while the $q_j$ are isobaric of weight $j-1$;
for instance $p_1=0$, $p_2=-K$, $p_3=-K'$, $p_4=K^2-K''$.

\begin{proposition}\label{prop:closed}
For every $k\ge0$,
\[
 G_k=\sum_{j=0}^{k}\binom kj\,p_j\,G_0^{(k-j)},\qquad G_0=NK=\frac Ky .
\]
\end{proposition}
\begin{proof}
Put $\tau_k:=S_k/y$. Since $\kappa=-y'/y$ and $N=1/y$, Lemma
\ref{lem:recursion} gives
\[
 G_k=\frac{K^{(k)}}{y}-y'\tau_k ,\qquad
 \tau_{k+1}=\tau_k'+\frac{K^{(k)}}{y^2},\qquad \tau_0=0 .
\]
For $u_1,u\in I$ put $H(u_1,u):=K(u)\int_{u_1}^{u}y^{-2}$, and
$\hat\tau_k(u_1):=\partial_u^kH(u_1,u)\big|_{u=u_1}$. Then $\hat\tau_0=0$, and
since $\partial_{u_1}H(u_1,u)=-K(u)/y(u_1)^2$,
\[
 \hat\tau_k'(u_1)=\big(\partial_{u_1}\partial_u^kH+\partial_u^{k+1}H\big)\big|_{u=u_1}
 =-\frac{K^{(k)}(u_1)}{y(u_1)^2}+\hat\tau_{k+1}(u_1).
\]
Thus the $\hat\tau_k$ satisfy the same recursion as the $\tau_k$, and
$\hat\tau_k=\tau_k$. Consequently, for each fixed $u_1\in I$,
\[
 G_k(u_1)=\frac{d^k}{du^k}\Big|_{u=u_1}K(u)\Big[\frac1{y(u_1)}-y'(u_1)\int_{u_1}^{u}y^{-2}\Big].
\]
Let $Y$ be the solution of $Y''+KY=0$ on $I$ with $Y(u_1)=1$, $Y'(u_1)=0$. Since
$Y$ and $y$ solve the same equation, their Wronskian $Yy'-Y'y$ is constant,
equal to $y'(u_1)$; hence $(Y/y)'=-y'(u_1)/y^2$ and
\[
 \frac{Y(u)}{y(u)}=\frac1{y(u_1)}-y'(u_1)\int_{u_1}^{u}y^{-2}.
\]
As $K/y=G_0$, this gives
$G_k(u_1)=\dfrac{d^k}{du^k}\Big|_{u=u_1}G_0(u)\,Y(u)$. Finally,
$Y^{(j)}=p_jY+q_jY'$ gives $Y^{(j)}(u_1)=p_j(u_1)$, and the Leibniz rule yields
\[
 G_k(u_1)=\sum_{j=0}^k\binom kj\,p_j(u_1)\,G_0^{(k-j)}(u_1).
\]
Since $u_1\in I$ is arbitrary, the claim follows.
\end{proof}

\section{Proof of the main results}\label{sec:proof}
Throughout this section $(M,g)$ is as in the Introduction, and $\U$ is the open
set defined before Proposition \ref{prop:reduction}, now for any $r\ge1$.
We first isolate the algebraic core of the argument.

\begin{lemma}\label{lem:residue}
Let $P$ and $A$ be polynomials and $\KK:=A/P$. Let $u_0\in\mathbb C$ be a root of
$P$ of multiplicity $m\ge2$ such that, with $s=u-u_0$,
\[
 P=\pi_0s^m\big(1+O(s)\big),\qquad
 \KK=\kappa_0s^{-2}+O(s^{-1}),\qquad \kappa_0=-(m+2)(m+3),\quad \pi_0\neq0 .
\]
Let $p_j[\KK]$ denote the differential polynomials $p_j$ of Section
\ref{sec:chain} evaluated at $\KK$ in place of $K$, and define
$G_k:=\sum_{j=0}^k\binom kj\,p_j[\KK]\,P^{(k-j)}$. Then the rational function
$G_{m+1}$ has at $u_0$ a simple pole with residue
\[
 \pi_0\,\frac{(m+1)!\,(m+2)}{2m+5}\,
 \Big[\binom{2m+3}{m+1}+(-1)^{m+1}(m+3)\Big]\neq0 .
\]
\end{lemma}

The proof is a local computation with the model $\KK_\ast=\kappa_0s^{-2}$; it
is given in Appendix~\ref{app:residue}.

\begin{proposition}\label{prop:Czero}
$C\equiv0$ on $\U$. Consequently $\LL$ and $\N$ are parallel on $\U$, and
$\U$ is locally isometric to a Riemannian product $\Sigma^2\times\mathbb
R^{n-2}$, where $\Sigma^2$ is a surface with Gaussian curvature $K$.
\end{proposition}

\begin{proof}
\textit{Case $r\ge3$.} Suppose, on the contrary, that the open set
$\V=\{p\in\U:\ C_p\neq0\}$ were non-empty, and fix an integral curve of $e_1$ in
$\V$ as in Section \ref{sec:chain}.

First we prove that, in that case,  $K$ and $y$ are rational.
By Lemma \ref{lem:words}, $P:=G_0=NK$ is a nonzero polynomial in $u$ of degree
$d\le r-2$. Moreover, $G_2$ is a polynomial because $r\ge3$. By Proposition
\ref{prop:closed}, $G_2=P''-KP$. Hence $A:=KP=P''-G_2$ is a nonzero polynomial,
and
\[
 K=\frac AP,\qquad y=\frac KP=\frac A{P^2}.
\]
We regard these expressions as rational functions of a complex variable $u$.
Since $y''+Ky=0$ holds on a real interval, it holds identically as an identity
of rational functions. Furthermore, Proposition \ref{prop:closed} shows that
every $G_k$ is rational, whereas Lemma \ref{lem:words} implies that, for
$0\le k\le r-1$, it coincides with a polynomial. We know ascertain the degree of the polynomial $A$.
Multiplying $y''+Ky=0$ by $P^4$ we get
\begin{equation}\label{eq:poly}
 A''P^2-4A'PP'-2APP''+6AP'^2+A^2P=0 .
\end{equation}
Let $\delta=\deg A$, with leading coefficients $\alpha$ of $A$ and $\pi$ of
$P$. The last term has degree $2\delta+d$ and leading coefficient
$\alpha^2\pi\neq0$. The sum of the first four has degree at most $\delta+2d-2$,
with coefficient $\alpha\pi^2(2d-\delta)(2d-\delta+1)$ in degree $\delta+2d-2$.
If $\delta>d-2$, then $2\delta+d>\delta+2d-2$, and the leading term of $A^2P$
cannot cancel. Hence $\delta\le d-2$, which in particular implies $d\ge2$.
Consequently $(2d-\delta)(2d-\delta+1)>0$, so the first four terms have exact
degree $\delta+2d-2$. If $\delta<d-2$, this degree is strictly greater than
$2\delta+d$, and their leading term cannot cancel either. Thus $\delta=d-2$, that is, $\deg A=\deg P-2$.

The next step is to prove that $K$ has a double pole at the roots of $P$, in order to be able to apply Lemma \ref{lem:residue} . Let $u_0\in\mathbb C$ be a root of $P$ of multiplicity $m\geq 1$. Since $\deg A<\deg P$, the order $\nu$ of $u_0$ as a root of $A$ can only be a natural number between $0$ and $m-1$ (otherwise $P$ would divide $A$). Then $y=y_0(u-u_0)^{\varepsilon}(1+O(u-u_0))$ with
$\varepsilon=\nu-2m=(\nu -m) -m \le-2$, and $K=\kappa_0(u-u_0)^{\nu-m}(1+O(s))$ with
$\nu-m\le-1$ and $y_0\kappa_0\neq0$. In $y''+Ky=0$ the leading term is
$$\varepsilon(\varepsilon-1)y_0(u-u_0)^{\varepsilon-2}+\kappa_0y_0(u-u_0)^{\varepsilon+\nu-m}$$
with $\varepsilon(\varepsilon-1)\neq0$.
This must vanish, which forces
\[
 \nu-m=-2,\qquad \varepsilon=-(m+2),\qquad \kappa_0=-(m+2)(m+3).
\]
Since $\nu\ge0$, this implies $m\ge2$. Hence, we are in the situation of Lemma \ref{lem:residue} with $\KK=A/P=K$ and by Proposition \ref{prop:closed} the functions $G_k$ of Section \ref{sec:chain} are those of that lemma.
 Hence $G_{m+1}$ has a nonzero residue at $u_0$. But $m\le d\le r-2$,
so $m+1\le r-1$ and $G_{m+1}$ is a polynomial by Lemma \ref{lem:words}. 
This contradiction shows that $\V=\emptyset$, i.e.\ $C\equiv0$ on $\U$.

\medskip\noindent\textit{Case $r=2$.} This is Lemma \ref{lem:r2}.

\medskip\noindent\textit{Case $r=1$.} Write $t=t^{(1)}$. Differentiating
$\nab\RR=t\otimes\RR$ gives
\[
 \nab^2\RR=t^{(2)}\otimes\RR,\qquad t^{(2)}:=\nab t+t\otimes t,
\]
so $(M,g)$ satisfies \eqref{eq:def} with $r=2$. Let
$\U_2:=\{p:\ t^{(2)}_p\neq0,\ \RR_p\neq0\}$. By the case $r=2$, the curvature
operator has rank one on $\U_2$ and $C\equiv0$ there. We claim that $\U_2$ is dense in $\U$. Otherwise, there
would exist a non-empty open subset $\mathcal O\subset\U$
disjoint from $\U_2$. Since $R\neq0$ on $\U$, this implies
$t^{(2)}=0$ on $\mathcal O$, and hence $\nabla^2R=0$ there.
By \cite{tanno}, this gives $\nabla R=0$ on $\mathcal O$,
contradicting $\mathcal O\subset\U$.
Consequently the curvature operator has rank at most one on $\U$, hence rank
exactly one, since $\RR\neq0$ there. By the algebraic Bianchi identity a rank
one curvature operator has the form \eqref{eq:typeII} with $F$ simple, so $\LL$,
$\N$ and $C$ are defined and continuous on $\U$. Since $\U\cap\U_2$ is dense in
$\U$ and $C=0$ there, $C\equiv0$ on $\U$.

\medskip
Finally, $C\equiv0$ means that $\N$, and hence $\LL$, is parallel on $\U$. By
Lemma \ref{lem:parallelL} the leaves of $\N$ are flat, and the local de Rham
theorem gives the stated product structure, the factor tangent to $\LL$ having
Gaussian curvature $K$.
\end{proof}

\begin{theorem}\label{thm:main}
Let $(M,g)$ be a connected Riemannian manifold of dimension $n\ge3$ whose local
holonomy group acts irreducibly at every point. If $(M,g)$ is $r$-th order
recurrent for some $r\ge1$, then it is locally symmetric.
\end{theorem}

\begin{proof}
If $\U=\emptyset$, then $t^{(r)}\otimes\RR\equiv0$, i.e.\ $\nab^r\RR=0$, and
$(M,g)$ is locally symmetric by \cite{lichnerowicz,tanno} (see Theorem 6(a) of
\cite{senovilla}). If $\U\neq\emptyset$, Proposition \ref{prop:Czero} shows that
$\LL$ is a parallel distribution of rank two on the open set $\U$. Since
$n\ge3$, it is a proper non-trivial parallel distribution, which contradicts the
irreducibility of the local holonomy.
\end{proof}

\begin{theorem}\label{thm:structure}
Let $(M,g)$ be a connected Riemannian manifold of dimension $n\ge3$ which is
$r$-th order recurrent, $\nab^{r}\RR=t\otimes\RR$, for some $r\ge1$. Then
$(M,g)$ is either locally symmetric or locally a flat extension of a surface.
More precisely, if $(M,g)$ is not locally symmetric, then:
\begin{enumerate}
\item[(a)] on $\operatorname{int}\{t=0\}$ one has $\nab\RR=0$;
\item[(b)] the curvature operator has rank at most one on all of $M$;
\item[(c)] on the open set $\{\RR\neq0\}$, $(M,g)$ is locally isometric to a
  product $\Sigma^2\times\mathbb{R}^{n-2}$. On $\{t\neq0\}\cap\{\RR\neq0\}$ the
  curvature $K$ of $\Sigma^2$ satisfies $\nab^{r}K\neq0$, whereas on the
  components of $\operatorname{int}\{t=0\}$ where $\RR\neq0$, $K$ is constant.
\end{enumerate}
\end{theorem}

\begin{proof}
If $(M,g)$ is locally symmetric there is nothing to prove, so assume it is not.
Let $\U$ be as above and $\Omega:=M\setminus\overline{\U}$. At every point of
$M\setminus\U$ we have $t\otimes\RR=0$, so $\nab^r\RR=0$ on the open set
$\Omega$, and $\nab\RR=0$ on $\Omega$ by \cite{tanno}. Since
$\operatorname{int}\{t=0\}\subset\Omega$, this proves (a).

If $\U=\emptyset$, then $\Omega=M$ and $\nab\RR=0$ on $M$, contrary to the
assumption that $(M,g)$ is not locally symmetric. Hence $\U\neq\emptyset$. (Conversely, if $\U\neq\emptyset$ then
$\nab^r\RR=t\otimes\RR\not\equiv0$, so $(M,g)$ is not locally symmetric: the
hypothesis is equivalent to $t\otimes\RR\not\equiv0$.)

(b) On $\U$ the curvature operator has rank one, by Proposition
\ref{prop:reduction} for $r\ge2$ and by the proof of Proposition
\ref{prop:Czero} for $r=1$; hence its rank is at most one on $\overline{\U}$.
Let $\Omega_0$ be a connected component of $\Omega$. Since $M$ is connected and
$\emptyset\neq\U\subset M\setminus\Omega_0$, the frontier of $\Omega_0$ is
non-empty, and it is contained in $\partial\U$, because $\Omega$ is open and
$\U$ is an open set disjoint from $\overline{\Omega}$. As $\nab\RR=0$ on
$\Omega_0$, parallel transport shows that the eigenvalues of the curvature
operator are constant on $\Omega_0$; by continuity they coincide with those at
a frontier point, where the rank is at most one. Hence the rank is at most one
on $\Omega$, and therefore on $M=\U\cup\partial\U\cup\Omega$.

(c) On the open set $M':=\{\RR\neq0\}$ the rank is exactly one, so, by the
algebraic Bianchi identity, \eqref{eq:typeII} holds locally on $M'$ with $F$
simple, and $\LL$, $\N$ and $C$ are defined and continuous on $M'$. We have
$C=0$ on $\U$ by Proposition \ref{prop:Czero}, and on $\Omega\cap M'$ because
there $\nab\RR=0$ forces the image of the curvature operator, spanned by $F$, to
be parallel, so that $\nab F=0$ and $C=0$ by Lemma \ref{lem:gradF}. Since
$\partial\U$ has empty interior, $(\U\cup\Omega)\cap M'$ is dense in $M'$, hence
$C\equiv0$ on $M'$. As in the proof of Proposition \ref{prop:Czero}, $M'$ is
locally isometric to $\Sigma^2\times\mathbb R^{n-2}$, and $F$ is parallel. Then
$\nab^r\RR=(\nab^rK)\,F\otimes F$, and comparing with $t\otimes\RR=tK\,F\otimes F$
gives $\nab^rK=tK$, which does not vanish on $\{t\neq0\}\cap\{\RR\neq0\}$.
Finally, on a component of $\operatorname{int}\{t=0\}$ where $\RR\neq0$ we
have $\nab\RR=0$, so $K=\RR(e_1,e_2,e_1,e_2)$ is locally constant.

Finally, (b) and (c) say precisely that $(M,g)$ is locally a flat extension of
a surface, which proves the first statement.
\end{proof}

\begin{remark}[Relation between Theorems \ref{thm:main} and \ref{thm:structure}]\label{rem:relation}
By de Rham's theorem, a Riemannian manifold is locally a product of a flat
factor and irreducible factors. Theorem \ref{thm:structure} is consistent with
the following reading: an irreducible factor of dimension at least three on which
$\RR\neq0$ is itself $r$-th order recurrent and hence, by Theorem
\ref{thm:main}, locally symmetric, so that the only non-symmetric factors are
surfaces, and the rank condition (b) shows that at most one such factor can
occur (cf.\ Theorem 1 of \cite{senovilla}). The proof given here does not use the decomposition. It is direct, and it
also controls the points where $\RR=0$, at which the local decompositions of
neighbouring regions have to be matched.
\end{remark}

\begin{remark}[Sharpness]\label{rem:sharp}
The exceptional case in Theorem \ref{thm:structure} does occur, and not only for
$r=1$. Let $\Sigma$ be the surface $du^2+G^2(u)dv^2$, $|u|<\epsilon$, where $G$
solves
\[
 G''=-uG,\qquad G(0)=1,\qquad G'(0)=0,
\]
so that $G>0$ for $\epsilon$ small and the Gaussian curvature of $\Sigma$ is
$K=-G''/G=u$. On $M=\Sigma\times\mathbb R^{n-2}$ the area form $F$ of $\Sigma$ is
parallel, hence $\nab^k\RR=(\nab^kK)\,F\otimes F$ for all $k$. Since
$\nab^2K=GG'\,dv\otimes dv$ and $G'(u)=-u^2/2+O(u^5)$,
\[
 \nab^2\RR=t^{(2)}\otimes\RR,\qquad t^{(2)}=\frac{GG'}{u}\,dv\otimes dv,
\]
where $t^{(2)}$ is smooth. On the other hand, at $u=0$ one has $\RR=0$ but
$\nab\RR=du\otimes F\otimes F\neq0$, so $M$ satisfies no condition of the form
$\nab\RR=t\otimes\RR$: it is proper of order two in the sense of
\cite{senovilla}, and it is not locally symmetric. In particular, the statement
about proper second-order recurrent Riemannian manifolds in Theorem 6(b) of
\cite{senovilla} requires, as for $r=1$, that flat extensions of surfaces be
excluded; Theorem \ref{thm:structure} provides the corrected statement for every
$r$.
\end{remark}

\begin{remark}
Although Lemma~\ref{lem:words} provides a whole family of polynomiality
conditions, the proof for $r\ge3$ uses only three of them, those corresponding
to $k=0$, $k=2$ and $k=m+1$. The mechanism is the following: the recurrence
conditions imply that $G_0=NK$ is a nonzero polynomial and, using the condition
for $k=2$, that $K$ is rational. The Jacobi equation $y''+Ky=0$ then forces $K$
to have, at some root $u_0$ of $G_0$ of multiplicity $m$, a double pole with
the specific leading coefficient $-(m+2)(m+3)$, that is,
\[
 K=-\frac{(m+2)(m+3)}{(u-u_0)^2}+O\big((u-u_0)^{-1}\big).
\]
For $k\le m$ the functions $G_k$ are regular at $u_0$, so this pole produces no
contradiction with their polynomiality; the first function in which it shows up
is $G_{m+1}$, which by Lemma~\ref{lem:residue} has a simple pole at $u_0$,
contradicting Lemma~\ref{lem:words}.
\end{remark}

\begin{remark}
The argument is local. It uses neither completeness, nor the Frenet-type
completion of the frame (Remark~\ref{rem:frenet}), nor the normal form of the
metric obtained in \cite{VanHookThesis,VanHook}. Its geometric input consists of
the curvature decomposition of Proposition~\ref{prop:reduction}, the nilpotency
of the splitting tensor (Lemma~\ref{lem:nilpotent}), the first-order relations
of Lemmas~\ref{lem:structure} and~\ref{lem:gradF}, valid on the whole open set
$\V$, and the recurrence condition.
\end{remark}

\section*{Acknowledgments}
In the preparation of this work I used Claude (Anthropic) to check some calculations (especially the Lemma in the Appendix) 
and to  typeset parts of the text. I have devised and checked all statements, and I take full
responsibility for the content.

\appendix

\section{Proof of Lemma \ref{lem:residue}}\label{app:residue}

\begin{proof}[Proof of Lemma~\ref{lem:residue}]
Since $p_j$ is isobaric of weight $j$ and $\KK$ has a double pole at $u_0$,
\[
 p_j[\KK]=\mathsf P_js^{-j}+O(s^{1-j}),
\]
where $\mathsf P_js^{-j}$ is the value of $p_j$ on the model
$\KK_\ast=\kappa_0s^{-2}$. Moreover, taking the $i$-th derivative, $P^{(i)}=\pi_0\frac{m!}{(m-i)!}s^{m-i}+O(s^{m-i+1})$
for $i\le m$, while $P^{(i)}=O(1)$ for $i>m$. Hence every term
$\binom{m+1}j p_j[\KK]P^{(m+1-j)}$ is $O(s^{-1})$, its coefficient of $s^{-1}$ only
depends on $\kappa_0$ and $\pi_0$, and for $j=0$ it is $P^{(m+1)}=O(1)$, which
contributes nothing. Consequently the residue of $G_{m+1}$ is the one computed
for the model $\KK_\ast=\kappa_0s^{-2}$, $P_\ast=\pi_0s^m$. For the model, the solutions of $Y''+\KK_\ast Y=0$ are $u^{m+3}$ and
$u^{-(m+2)}$, since $u_0=0$ and $s=u$. Fix $u_1\neq0$. The solution with
$Y(u_1)=1$, $Y'(u_1)=0$ is
\[
 Y(u)=\frac{(m+2)(u/u_1)^{m+3}+(m+3)(u/u_1)^{-(m+2)}}{2m+5},
\]
and, writing $u=u_1(1+x)$,
\[
 P_\ast(u)Y(u)=\frac{\pi_0u_1^m}{2m+5}
 \Big[(m+2)(1+x)^{2m+3}+(m+3)(1+x)^{-2}\Big].
\]
The coefficient of $x^k$ is
$\beta_k=(m+2)\binom{2m+3}{k}+(-1)^k(m+3)(k+1)$. Since $Y^{(j)}(u_1)=p_j[\KK_\ast](u_1)$,
the Leibniz rule gives
\[
 G_k(u_1)=\frac{d^k}{du^k}\Big|_{u=u_1}P_\ast(u)\,Y(u)
 =\frac{k!\,\beta_k}{u_1^k}\cdot\frac{\pi_0u_1^m}{2m+5},
\]
because $d/du=u_1^{-1}\,d/dx$. Thus the model gives
$G_k(u_1)=\pi_0\,k!\,\beta_k\,u_1^{m-k}/(2m+5)$ for every $u_1\neq0$. For
$k=m+1$ this is the stated residue, which does not vanish because
$\binom{2m+3}{m+1}\ge\binom{2m+3}{2}=(2m+3)(m+1)>m+3$.
\end{proof}

\end{document}